\documentclass[11pt]{amsart}
\usepackage{amsfonts,amssymb,amscd, graphicx}
\usepackage[all]{xy}

\newcommand{\LP}{\operatornamewithlimits{\overrightarrow{\textstyle\prod}}}
\newcommand{\RP}{\operatornamewithlimits{\overleftarrow{\textstyle\prod}}}

\newcommand{\IR}{\mathbb R}
\newcommand{\IN}{\mathbb N}
\newcommand{\IS}{\mathbb S}
\newcommand{\IZ}{\mathbb Z}
\newcommand{\ID}{\mathbb D}
\newcommand{\IK}{\mathbb K}
\newcommand{\IP}{\mathbb P}
\newcommand{\II}{\mathbb I}
\newcommand{\IT}{\mathbb T}
\newcommand{\IM}{\mathbb M}

\newcommand{\IA}{\mathbb A}
\newcommand{\w}{\omega}
\newcommand{\U}{\mathcal U}
\newcommand{\V}{\mathcal V}

\newcommand{\M}{\mathcal M}

\newcommand{\Ra}{\Rightarrow}
\newcommand{\supp}{\operatorname{supp}}

\newcommand{\cov}{\mathrm{cov}}
\newcommand{\St}{\mathcal St}
\newcommand{\ulim}{\mathrm u\mbox{-}\kern-2pt\varinjlim}
\newcommand{\glim}{\mathrm g\mbox{-}\kern-2pt\varinjlim}

\newcommand{\Int}{\mathrm{int}}

\newcommand{\HH}{\mathcal H}

\newcommand{\EE}{\mathcal E}

\newcommand{\id}{\mathrm{id}}

\newcommand{\cl}{\mathrm{cl}}
\newcommand{\bd}{\mathrm{bd}}
\newcommand{\tint}{\mathrm{int}}

\newtheorem{theorem}{Theorem}
\newtheorem{lemma}{Lemma}
\newtheorem*{claim}{Claim}

\theoremstyle{definition}
\newtheorem{definition}{Definition}
\newtheorem{remark}{Remark}

\title[On homeomorphism groups of non-compact surfaces]{On homeomorphism groups of non-compact surfaces, \\
endowed with the Whitney topology}

\author[T. Banakh]{Taras Banakh}
\address[T. Banakh]{Instytut Matematyki,
Jan Kochanowski University, Kielce, Poland, and\newline
 Department of Mathematics,
 Ivan Franko National University of Lviv, 79000, Ukraine}
 \email{tbanakh@yahoo.com}

\author[K. Mine]{Kotaro Mine}
 \address[K. Mine]{Institute of Mathematics,
 University of Tsukuba, Tsukuba, 305-8571, Japan}
 \email{pen@math.tsukuba.ac.jp}

\author[K. Sakai]{Katsuro Sakai}
 \address[K. Sakai]{Institute of Mathematics,
 University of Tsukuba, Tsukuba, 305-8571, Japan}
 \email{sakaiktr@sakura.cc.tsukuba.ac.jp}

\author[T. Yagasaki]{Tatsuhiko Yagasaki}
\address[T. Yagasaki]{Division of Mathematics,
 Kyoto Institute of Technology, Kyoto, 606-8585, Japan}
\email{yagasaki@kit.ac.jp}

\thanks{The last author was supported by JSPS Grant-in-Aid for Scientific Research (No.22540081).}


\begin{document}
\begin{abstract} We prove that for any non-compact connected surface $M$ the group $\HH_c(M)$ of compactly suported homeomorphisms of $M$ endowed with the Whitney topology is homeomorphic to $\IR^\infty\times l_2$ or $\IZ\times\IR^\infty\times l_2$.
\end{abstract}
\maketitle

\section{Introduction}

This paper is devoted to studying the topological structure of (the identity component $\HH_0(M)$ of) the homeomorphism group $\HH(M)$ of a non-compact connected surface $M$.
By a {\em surface} we understand a $\sigma$-compact 2-manifold $M$ possibly with boundary $\partial M$.
By \cite{Moise}, each surface $M$ admits a combinatorial triangulation, unique up to PL-homeomorphisms.
If necessary, we fix a triangulation of $M$ and regard $M$ as a PL $2$-manifold.
A {\em subpolyhedron} of $M$ means a subpolyhedron
 with respect to a PL-structure on $M$.

For a surface $M$ let $\HH(M)$ denote the homeomorphism group of $M$ endowed with the Whitney topology.
This topology is generated by the base consisting of the sets
$$\Gamma_U=\{f\in\HH(M):\Gamma_f\subset U\},$$
where $U$ runs over open subsets of $M\times M$ and $\Gamma_f=\{(x,f(x)):x\in M\}$ stands for the graph of a function
$f:M\to M$. If $M$ is compact, then the Whitney topology on $\HH(M)$ coincides with the compact-open topology.  By \cite[Proposition 3.2]{BMSY1}, $\HH(M)$ is a topological group.

Given a subset $K\subset M$, consider the subgroup $$\HH(M;K)=\{f\in\HH(M):f|_K=\id_K\}\subset \HH(M).$$
Let $\HH_0(M;K)$ be the identity connected component of $\HH(M;K)$.
This is the largest connected subset that contains the neutral element $\id_M$ of the group $\HH(M;K)$.
By \cite[Proposition 3.3]{BMSY1} $\HH_0(M;K)$ lies in the subgroup
$\HH_c(M;K)$ of $\HH(M;K)$ that consists of all homeomorphisms $f\in \HH(M;K)$ having compact support $$\supp(f)=\cl_M\{x\in M:f(x)\ne x\}\subset M.$$
We write $\HH(M)$, $\HH_0(M)$ and $\HH_c(M)$ instead of $\HH(M;\emptyset)$, $\HH_0(M;\emptyset)$, and $\HH_c(M;\emptyset)$.

The local topological structure of the groups $\HH_0(M;K)$ and $\HH_c(M;K)$ was studied in \cite{BMSY1} for the case when $K$ is a subpolyhedron in a surface $M$. It was shown in \cite{BMSY1} that the topological group $\HH_c(M;K)$ is
\begin{itemize}
\item an $l_2$-manifold if $\cl_M(M\setminus K)$ is compact;
\item an $(\IR^\infty\times l_2)$-manifold, otherwise.
\end{itemize}
Here $\IR^\infty$ is the direct limit of the tower
$$\IR\subset\IR^2\subset\IR^3\subset\cdots$$ of Euclidean spaces, where each $\IR^n$ is identified with the hyperplane $\IR^n\times\{0\}$ in $\IR^{n+1}$.

Another result of \cite{BMSY1} says that the topological group $\HH_c(M;K)$ is locally contractible and hence the connected component $\HH_0(M;K)$ is an open subgroup of $\HH_c(M;K)$. Consequently, $\HH_c(M;K)$ is homeomorphic to the product $\HH_0(M;K)\times\M_c(M;K)$ of the connected group $\HH_0(M;K)$ and the discrete quotient group
$$\M_c(M;K)=\HH_c(M;K)/\HH_0(M;K),$$
which can be called the {\em mapping class group} of the pair $(M, K)$.
Therefore, the topological type of the group $\HH_c(M;K)$ is completely determined by that of $\HH_0(M;K)$ and the cardinality of $\M_c(M;K)$.

The topological structure of the group $\HH_0(M;K)$ is well-understood in case of a compact connected surface $M$. The following classification result belongs to M.~Hamstrom \cite{Ham} (cf.\ \cite{Yag2}).

\begin{theorem}[Hamstrom]\label{Ham}
For a subpolyhedron $K\subsetneqq M$ in a compact connected surface $M$ the group $\HH_0(M;K)$ is homeomorphic to
\begin{itemize}
\item $\IS^1\times l_2$  if the pair $(M,K)$ is homeomorphic to one of the pairs: \smallskip
\newline \hspace*{20mm}
$(\ID,\emptyset)$, $(\ID,\{0\})$, $(\IA,\emptyset)$, $(\IS^2,1pt)$,
$(\IS^2,2pt)$, $(\IP,1pt)$, $(\IM,\emptyset)$ or $(\IK,\emptyset)$;
\smallskip
\item $\IT\times l_2$  if $(M,K)$ is homeomorphic to $(\IT,\emptyset)$;
\item $\mathit{SO}(3)\times l_2$ if $(M,K)$ is homeomorphic to $(\IS^2,\emptyset)$
 or $(\IP,\emptyset)$;
\item $l_2$ in all other cases.
\end{itemize}
\end{theorem}

Here $\ID$ denotes the disk, $\IA=\ID\#\ID$ the annulus, $\IP$ the projective plane, $\IM=\IP\#\ID$ the M\"obius band, $\IK=\IP\#\IP$ the Klein bottle,  and $\IT$ the torus.
As expected, $\IS^n$ denotes the $n$-dimensional sphere. For two surfaces $M,N$ by $M\#N$ we denote their connected sum. The connected sum $M\#\ID$ is homeomorphic to $M$ with removed open disk.

The classification Theorem~\ref{Ham} is completed by the following theorem that will be  proved in Section~\ref{spH0}.

\begin{theorem}\label{H0}
For a subpolyhedron $K\subsetneqq M$ in a non-compact connected surface $M$ the group $\HH_0(M;K)$ is homeomorphic to
\begin{enumerate}
\item
$l_2$ if $\cl_M(M\setminus K)$ is compact;
\item
$\IR^\infty\times l_2$, otherwise.
\end{enumerate}
\end{theorem}

For non-compact graphs a counterpart of Theorem~\ref{H0} was proved in \cite{BMS}.

Next, we calculate the cardinality of the mapping class group $\M_c(M)=\HH_c(M)/\HH_0(M)$ of a connected surface $M$. For compact connected surfaces the following classification is known, see \cite{Bir}, \cite[Ch.2]{Farb}, \cite{Kork}, \cite{Li1}, \cite{Li2}, \cite{Szep}, \cite[\S4.4]{VMP}.

\begin{theorem}\label{MCG-c} For a compact connected surface $M$ the mapping class group $\M_c(M)$ is isomorphic to
\begin{itemize}
\item the trivial group if $M$ is homeomorphic to $\IP$;
\item $\IZ_2$ if $M$ is homeomorphic to $\IS^2$, $\ID$, or $\IM$;
\item $D_4$ if $M$ is homeomorphic to $\IA$ or $\IK$;
\item $D_8$ if $M$ is homeomorphic to $\IM\#\ID$, the M\"obius band with a hole;
\item $D_{12}$ if $M$ is homeomorphic to $\ID\#\ID\#\ID$, the disk with two holes;
\item $GL(2,\IZ)$ if $M$ is homeomorphic to the torus $\IT$;
\item an infinite group in all other cases.
\end{itemize}
\end{theorem}

Here $D_{2n}=\langle a,b\mid a^n=b^2=baba=1\rangle$ is the dihedral group (of isometries of the regular $n$-gon) having $2n$ elements. It should be mentioned that the group $D_2$ is isomorphic to ($\cong$) the cyclic group $\IZ_2=\IZ/2\IZ$, $D_4\cong\IZ_2\times\IZ_2$, $D_6$ is isomorphic to the symmetric group $\Sigma_3$, and $D_{12}\cong D_6\times\IZ_2\cong \Sigma_3\times\IZ_2$.

We complete Theorem~\ref{MCG-c} by the following theorem that will be proved in Section~\ref{sMCG}.

\begin{theorem}\label{MCG}
For a non-compact connected surface $M$,
 the following conditions are equivalent:
\begin{enumerate}
\item $\M_c(M)$ is trivial;
\item $\M_c(M)$ is a torsion group\footnote{A group $G$ is called {\em torsion} if each element of $G$ has finite order.};
\item $M$ is homeomorphic to $X\setminus K$,
 where $X = \IA$, $\ID$ or $\IM$,
 and $K$ is a non-empty compact subset of a boundary circle of $X$.
\end{enumerate}
\end{theorem}

To distinguish the surfaces appearing in Theorems~\ref{Ham}, \ref{MCG-c}, \ref{MCG}, let us introduce the following terminology:

\begin{definition} A surface $M$ is called {\em exceptional} if it is homeomorphic to one of the surfaces: $\IS^2$, $\IT$, $\IP$, $\IK$, $\IM\#\ID$, $\ID\#\ID\#\ID$ or
$X\setminus K$, where $X=\IA$, $\ID$ or $\IM$ and $K$ is a (possibly empty) compact subset of a boundary circle of $X$.
\end{definition}

Unifying Theorems~\ref{Ham}--\ref{MCG}, we get the following table describing all possible topological types of the homeomorphism groups $\HH_0(M)$ and $\HH_c(M)$ of a connected surface $M$. In this table cardinal numbers (identified with the sets of smaller ordinals) are endowed with the discrete topology.

\smallskip
$$\begin{tabular}{|@{\quad}c@{\quad}|@{\quad}c@{\quad}|@{\quad}c@{\quad}|@{\quad}c@{\quad}|}
\hline
 & \\[-13.5pt]
$\phantom{I_{I_I}}M\phantom{I^{I^I}}$ & $\M_c(M)$ & $\HH_0(M)$ & $\HH_c(M)$ \\[-12pt]
 & \\
\hline
 & \\[-14pt]
$\phantom{I^{I^I}}\IP\phantom{I^{I^I}}$ & $\IZ_1$ & $SO(3)\times l_2$ & $SO(3)\times l_2$ \\
$\IS^2$ & $\IZ_2$ & $SO(3)\times l_2$ & $2\times SO(3)\times l_2$ \\
$\ID$, $\IM$ & $\IZ_2$&  $\IS^1\times l_2$ & $2\times \IS^1\times l_2$\\
$\IA$, $\IK$ & $D_4$ &$\IS^1\times l_2$ & $4\times \IS^1\times l_2$\\
$\IM\#\ID$ & $D_8$& $l_2$ & $8\times l_2$\\
$\ID\#\ID\#\ID$ & $D_{12}$ &$l_2$ & ${12}\times \IS^1\times l_2$\\
$\IT$ & $GL(2,\IZ)$ & $\IT\times l_2$ & $\IN\times \IT\times l_2$\\
compact, non-exceptional & infinite & $l_2$ & $\IN\times l_2$\\
non-compact, exceptional & trivial & $\IR^\infty\times l_2$ & $\IR^\infty\times l_2$ \\
non-compact, non-exceptional & infinite & $\IR^\infty\times l_2$ & $\IN\times \IR^\infty\times l_2$ \\[-13.5pt]
 & \\
\hline
\end{tabular}$$
\medskip

It is interesting to compare the classification Theorem~\ref{H0}  with the following result on the compact-open topology due to the last author \cite{Yag2}.

\begin{theorem}[Yagasaki] For a non-compact connected  surface $M$ the identity component  of the homeomorphism group $\HH(M)$ endowed with the compact-open topology is homeomorphic to:
\begin{enumerate}
\item
 $\IS^1\times l_2$ if $M$ is homeomorphic to $\IR^2$, $\IS^1\times\IR$,
 $\IS^1\times[0,\infty)$ or $\IM\setminus\partial \IM$;
\item
 $l_2$ in all other cases.
\end{enumerate}
\end{theorem}

\section{Recognizing topological groups homeomorphic to $\IR^\infty\times l_2$}

In this section we recall a criterion  for recognizing topological groups homeomorphic to $\IR^\infty\times l_2$, which is proved in \cite{BMRSY} and is based upon the topological characterization of the space $\IR^\infty\times l_2$ given in \cite{BR3}.
First we recall some necessary definitions.

A subgroup $H$ of a topological group $G$ is called {\em locally topologically complemented} (LTC) in $G$ if $H$ is closed in $G$ and the quotient map $q:G\to G/H=\{xH:x\in G\}$ is a locally trivial bundle, which happens if and only if $q$ has a local section at some point of $G/H$.
Here, a {\em local section} of a map $q: X \to Y$ at a point $y \in Y$ means a continuous map $s: U \to X$ defined on a neighborhood $U$ of $y$ in $Y$ such that $q\circ s=\id_U$. (As usual, $\id_U$ denotes the inclusion map $U \subset Y$.) We understand that the distinguished point of a group $G$ is its neutral element $e$, and for a subgroup $H$ the distinguished element of the quotient space $G/H$ is the coset $\bar e=e H$.

\begin{lemma}\label{l_ltc}
Suppose $G$ is a topological group and $K \subset H$ are closed subgroups of $G$.
\begin{itemize}
\item[(1)] If $G$ is metrizable, then so is the quotient space $G/H$.
\item[(2)] If $K$ is LTC in $H$ and $H$ is LTC in $G$, then $K$ is LTC in $G$.
\item[(3)] If $H$ is LTC in $G$, then the map $\pi :  G/K \to G/H$, $\pi(gK) = gH$,
 is a locally trivial bundle with the fiber $H/K$.
\end{itemize}
\end{lemma}

\begin{proof}
(1) Being metrizable, the group $G$ admits a right invariant metric $d$ generating the topology of $G$.
Then the topology of the quotient space $G/H$ is generated by
the metric $\rho$ defined by
$$\rho(xH,yH)=\inf\{d(a,b):a\in xH,\; b\in yH\} \hspace{5mm} (xH,yH\in G/H).$$

(2)
By the assumption,
the projection $G \to G/H$ has a local section
 $\sigma : (W, \overline{e}) \to (G, e)$ at $\overline{e} \in G/H$ and
the projection $H \to H/K$ also has a local section
 $\tau : U \to H$ at $\overline{e} \in H/K$.
 Consider the projection $\pi : G/K \to G/H$.
The map $\sigma$ determines the map
 $\sigma_0 : \pi^{-1}(W) \to H/K$, $\sigma_0(x) = \sigma(\pi(x))^{-1}x$.
Since $\sigma_0(\overline{e}) = \overline{e} \in U$,
there exists an open neighborhood $V$ of $\overline{e}$ in $\pi^{-1}(W)$
 such that $\sigma_0(V) \subset U$.
The required local section $s : V \to G$ for the projection $G \to G/K$ is defined by
 $s(x) = \sigma(\pi(x)) \tau(\sigma_0(x))$.

(3)
The projection $G \to G/H$ has a local section
 $\sigma : W \to G$ at any point $x_0 \in G/H$.
The associated trivialization
 $\phi : W \times H/K \approx \pi^{-1}(W)$ is defined by
 $\phi(x, y) = \sigma(x)y$.
\end{proof}

\begin{remark}\label{covering} We can also show the following statements.
Suppose  $H$ is a closed subgroup of $G$ and $K$ is an open normal subgroup of $H$.
\begin{enumerate}
\item
The natural map $\xi:G/K \to G/H$, $\xi(gK) = gH$,
 is a covering map with fiber $H/K$.
\item
 $H$ is locally topologically complemented in $G$
 if and only if so is $K$ in $G$.
\end{enumerate}
\end{remark}

Following \cite{BMRSY}, we say that a topological group $G$ carries the {\it strong topology} with respect to a tower of subgroups
$$G_0\subset G_1\subset G_2\subset\cdots$$
if $G=\bigcup_{n\in\w}G_n$ and for any
neighborhood $U_n$ of the neutral element $e$ in $G_n$, $n \in \w$,  the group product
$$\LP_{n\in\w}U_n=\bigcup_{n\in\w} U_0U_1\cdots U_n$$is a neighborhood of $e$ in $G$.

Let us recall that a closed subset $A$ of a topological space $X$ is called a ({\em strong}) {\em $Z$-set in} $X$ if for any open cover $\U$ of $X$ there is a continuous map $f:X\to X$ such that $f$ is $\U$-near to the identity $\id_X:X\to X$ and (the closure $\overline{f(X)}$ of) the set $f(X)$ does not intersect $A$.
A point $x$ of $X$ is called a ({\em strong}) {\em $Z$-point} if the singleton $\{x\}$ is a (strong) $Z$-set in $X$.
Every point in an infinite-dimensional Hilbert manifold is a strong $Z$-point.

The following theorem proved in \cite{BMRSY} will be our main instrument in the proof of Theorem~\ref{H0}.

\begin{theorem}\label{BMRSY} A non-metrizable topological group $G$ is homeomorphic to the space $\IR^\infty\times l_2$ if $G$ carries the strong topology with respect to  a tower of subgroups $(G_n)_{n\in\w}$ such that each group $G_n$ is homeomorphic to $l_2$, is locally topologically complemented in $G_{n+1}$, and each $Z$-point of the quotient space $G_{n+1}/G_n$ is a strong $Z$-point.
\end{theorem}

\section{Proof of Theorem~\ref{H0}}\label{spH0}

Suppose $M$ is a non-compact connected surface and
$K\subsetneqq M$ is a subpolyhedron of $M$ (with respect to some PL-structure $\tau$ of $M$). 
The boundary and interior of $M$ as a manifold are denoted by $\partial M$ and ${\rm Int}\,M$, while
the topological closure, interior and frontier of a subset $A$ in $M$
are denoted by the symbols $\cl_M A$, ${\rm int}_M A$ and ${\rm bd}_M A$ respectively.
%
%

Case (1) follows from the following lemma.

\begin{lemma}\label{lm_2}
If $\cl_M(M\setminus K)$ is compact, then $\HH_0(M;K) \approx l_2$ and
$\HH_0(M;K)$ is open in $\HH(M;K)$.
\end{lemma}

\begin{proof}
There exists a compact connected PL 2-submanifold $N$ of $(M, \tau)$ such that
$\cl_M(M\setminus K) \subset N$.
By the Hamstrom's Theorem~\ref{Ham}, $\HH_0(N; N\cap K) \approx l_2$.
In fact, $N \cap K$ contains the frontier ${\rm bd}_M N$,
which is a 1-manifold and is nonempty since $M$ is connected.
By \cite{Yag2} $\HH(N;N\cap K)$ is an $l_2$-manifold (in particular, it is locally connected).
The conclusion now follows from
$(\HH(M;K), \id_M) \approx (\HH(N;N\cap K), \id_N)$.
\end{proof}

Below we assume that $\cl_M(M\setminus K)$ is non-compact.
Then the surface $M$ can be represented as the countable union $M=\bigcup_{n\in\w}M_n$ of compact subpolyhedra $M_n$, $n\in\w$, of $(M, \tau)$ such that
$M_n \subset \Int_M M_{n+1}$ and $\Int_M M_n \not\subset K$ for each $n \in\w$.
Consider the subpolyhedra $\check M_n = M\setminus \Int_M M_n$ and $K_n = K \cup \check M_n$ $(n\in\w)$ of $(M, \tau)$.
Then we obtain the subgroups $H = \HH_c(M;K)$ and $H_n = \HH(M;K_n)$ $(n\in\w)$ of $\HH(M; K)$.
Since $K_n \supset K_{n+1} \supset K$ we have $H_n \subset H_{n+1}$ $(n\in\w)$ and  $H = \bigcup_{n\in\w} H_n$.

Our main concern is the  identity connected components $G = \HH_0(M;K)$ and $G_n = \HH_0(M;K_n)$ $(n\in\w)$.
Note that $G$ is not metrizable and it coincides with the identity connected component of the group $H$.
In Lemmas~\ref{lm_3}\,--\,\ref{lm_5} below we show that the tower of subgroups $G_n$, $n\in\w$, of the group $G$ satisfies the conditions in Theorem~\ref{BMRSY}.
This implies Case (2).

\begin{lemma}\label{lm_3}
{\rm (1)} $G = \bigcup_{n\in\w} G_n$.
{\rm (2)} $G_n \approx l_2$ and $G_n$ is open in $H_n = \HH(M;K_n)$ for each $n \in \w$.
\end{lemma}

\begin{proof}
(1) First note that $G = \HH_0(M;K)$ is path-connected since
the group $\HH_c(M;K)$ is locally path-connected
by \cite[Theorem 6.5]{BMSY1} and
$G$ is the identity connected component of $\HH_c(M;K)$.
Hence any $h \in G$ can be joined to $\id_M$ by an arc $A$ in $G$.
By \cite[Proposition 3.3]{BMSY1} the compact subset $A$ lies in $\HH(M;\check M_n)$ for some $n \in \w$.
Since
$$\id_M \in A \subset \HH(M;\check M_n)\cap \HH_0(M;K) \subset \HH(M;K_{n}),$$
the connectedness of $A$ implies that $h \in A \subset G_n = \HH_0(M;K_{n})$.
\smallskip

(2) The assertion follows from Lemma~\ref{lm_2} since
$K_n \subsetneqq M$ and $M \setminus K_n$ is included in the compact subset $M_n$.
\end{proof}

%
%

To check the LTC condition and the $Z$-point property in Theorem~\ref{BMRSY} for the tower $(G_n)_{n \in \w}$
we need a result on embedding spaces.
For closed subsets $L\subset N$ of $M$ let $\EE_L^*(N,M)$ denote the space of closed embeddings $f:N\to M$ such that $f|_L=\id_L$ and $f^{-1}(\partial M)=N\cap \partial M$. The space $\EE_L^*(N,M)$ is endowed with the compact open topology.
The following basic fact is due to Yagasaki \cite{Yag1}, and in its present form can be found in \cite[Theorem 6.7]{BMSY1}.
One should note that the assertion was verified in \cite{LM} in the important case that $K = \emptyset$ and
 $L$ is either a proper arc, an orientation-preserving circle or a compact 2-submanifold of $M$.

\begin{theorem}[Yagasaki]\label{Yag}
Suppose $L \subset N$ are two subpolyhedra of a surface $M$ with compact closure $\cl_M(N \setminus K)$.
\begin{enumerate}
\item
The restriction map
$$R :\HH(M,L) \to \mathcal E^*_L(N,M),
 \quad R(h) = h|_N$$
 has a local section $s : (\U, \id_N) \to (\HH_0(M,L),\id_M)$
 at the inclusion $\id_N : N \subset M$.
\item
\begin{itemize}
\item [(i)\,] ${\rm Im}\,R$ is an open neighborhood
 of the inclusion $\id_N$ in $\mathcal E^*_L(N,M)$.
\item[(ii)] The map $R :\HH(M,L) \to {\rm Im}\,R$
 is a principal $\HH(M,N)$-bundle. (In particular, the map $R$ is an open map.)
\end{itemize}
\item The space $\mathcal E^*_L(N,M)$
is an $l_2$-manifold if $\dim(N \setminus L) \ge 1$.
\end{enumerate}
\end{theorem}

\begin{lemma}\label{lm4} {\rm (1)} The subgroup $H_m$ is LTC in $H_n$ and the quotient space $H_n/H_m$ is an $l_2$-manifold for any $m \le n$ in $\w$.
\begin{itemize}
\item[(2)] The subgroup $G_m$ is LTC in $G_n$ and the quotient space $G_n/G_m$ is an $l_2$-manifold
(in particular, every point of $G_n/G_m$ is a strong $Z$-point) for any $m \le n$ in $\w$.
\end{itemize}
\end{lemma}

\begin{proof}
(1)
The group $H_n=\HH(M;K_n)$ acts continuously on the space $\EE^*_{K_n}(K_m,M)$ by the left composition.
Under this action
the subgroup $H_m$ is the isotropy group of $\id_{K_m}$
and the restriction map
$$R: H_n \to \EE^*_{K_n}(K_m,M),\;\;R:h\mapsto h|_{K_m},$$
coincides with the orbit map at $\id_{K_m}$.
Hence we have the factorization
$$\xymatrix@M+1pt{
 & H_n \ar[dl]_r \ar[dr]^R & \\
H_n/H_m \ar[rr]^{\phi}_{\cong} && {\rm Im}\,R & \hspace{-9mm} \subset \hspace{2mm} \EE^*_{K_n}(K_m,M),
}$$
where the map $r$ is the natural projection and the induced map $\phi$ is defined by
$\phi(hH) = h|_{K_m}$.
Note that the map $\phi$ is a homeomorphism since $\phi$ is bijective and both maps $r$ and $R$ are open continuous surjections.
Therefore, by Theorem~\ref{Yag} the map $r$ has a continuous section and $H_n/H_m$ is an $l_2$-manifold.

(2) First consider the subgroups $H_n \supset H_m \supset G_m$.
By Lemma~\ref{lm_3} the subgroup $G_m$ is open in $H_m$ and $H_m/G_m$ is discrete.
Hence $G_m$ is TC in $H_m$ and  by (1) and Lemma~\ref{l_ltc}\,(2) it follows that $G_m$ is LTC in $H_n$. From (1) and Lemma~\ref{l_ltc}\,(3) it also follows that the projection $\pi : H_n/G_m \to H_n/H_m$ is a locally trivial bundle with the discrete fiber $H_m/G_m$.
This implies that $H_n/G_m$ is also  an $l_2$-manifold. In fact, $H_n/G_m$ is locally homeomorphic to $l_2$ since $H_n/H_m$ is an $l_2$-manifold by (1) and
the metrizability of $H_n/G_m$ follows from that of $H_n$ and Lemma~\ref{l_ltc}\,(1).

Since $G_n$ is an open subgroup of $H_n$ and $G_m \subset G_n$, one sees that $G_n/G_m$ is an open neighborhood of
the distinguished element $\bar e=e G_m$ in $H_n/G_m$.
Hence any local section $s : (U, \bar e) \to (H_n, e)$ of the projection $H_n \to H_n/G_m$ restricts to
a local section of the projection $G_n \to G_n/G_m$. This means that $G_m$ is LTC in $G_n$.
Since $H_n/G_m$ is an $l_2$-manifold, so is the open subset $G_n/G_m$.
This completes the proof.
\end{proof}

It remains to check the strong topology condition for the tower $(G_n)_{n \in \w}$.
This would be derived from Theorem~\ref{Yag} and general arguments on transformation groups with strong topology in \cite[Section 5]{BMSY1}.
However, for the convenience  of the reader, here we include its direct self-contained proof and
avoid exceeding generality on related subjects.

Let us recall some notation. For covers $\U$, $\V$ of a space $M$ and a subset $A\subset M$ we put $\St(A,\V)=\cup\{V\in\V:A\cap V\ne\emptyset\}$,  $\St(\U,\V)=\{\St(U,\V):U\in\U\}$, and $\St(\U)=\St(\U,\U)$. We write $\U\prec\V$ if each set $U\in\U$ lies in some set $V\in\V$.

For two maps $f,g:X\to M$ we write $(f,g)\prec\U$ if for each point $x\in X$ the doubleton $\{f(x),g(x)\}$ lies in some subset $U\in\U$.

By \cite[Proposition 3.1]{BMSY1} the sets $B(\id_M,\U)=\{h\in\HH(M):(h,\id_M)\prec\U\}$, $\U\in\cov(M)$, form a neighborhood base of the Whitney topology at the neutral element $\id_M$ of the homeomorphism group $\HH(M)$. Here $\cov(M)$ denotes the family of all open covers of $M$.

\begin{lemma}\label{lm_5} For any neighborhood $U_n$ of $\id_M$ in $G_n=\HH_0(M;K_n)$, $n\in\w$,
the product $\LP_{n\in\w}U_n$ is a neighborhood of $\id_M$ in $\HH_c(M;K)$.
\end{lemma}

\begin{proof}
We may assume that each $U_n$ is symmetric (i.e., $U_n=U_n^{-1}$).
It suffices to show that the set
$$\Big(\LP_{n\in\w}U_n\Big)^{-1}=\RP_{n\in\w}U_n^{-1}=\RP_{n\in\w}U_n=\bigcup_{n \in w} U_n\cdots U_0$$
is a neighborhood of $\id_M$ in $\HH_c(M;K)$.

By Theorem~\ref{Yag}, for every $n\in\w$ the restriction map
$$R_n : G_n \to \EE^*_{K_n}(K_{n-1},M),\;\;R_n :h\mapsto h|_{K_{n-1}}$$
has a local section
$s_n : (V_n, \id_{K_{n-1}}) \to (G_n, \id_M)$ at $\id_{K_{n-1}}$.

Inductively we can find covers $\U_n, \V_n \in \cov(M)$, $n\in\w$, satisfying the next conditions for each $n\in\w$:
\begin{enumerate}
\item
\begin{itemize}
\item[(i)\,] $\St(\U_n)\prec\V_{n-1}$ (where $\V_{-1}=\{M\}$);
\item[(ii)] If $h \in \HH(M; K_n)$ and $(h, \id_M) \prec \U_n$, then $h \in U_n$;
\end{itemize}
\smallskip
\item
\begin{itemize}
\item[(i)\,] $\St(\V_n)\prec\U_n$;
\item[(ii)] If $f\in \EE^*_{K_n}(K_{n-1},M)$ and $(f,\id_{K_{n-1}})\prec \St(\V_n)$, then $f\in V_n$ and $(s_n(f),\id_M)\prec\U_n$.
\end{itemize}
\end{enumerate}

By the induction on $n \in \w$ we shall verify the following assertion:
\begin{itemize}
\item[$(\ast)_n$] If $h \in \HH(M; K_n)$ and $(h|_{K_{m-1}} , \id_{K_{m-1}}) \prec \St(\V_m,\St(\U_{n+1}))$ for each $m \in \w$ with $m \le n$, then
$h \in U_n \cdots U_0$.
\end{itemize}

$(\ast)_0$ : Since $h \in \HH(M; K_0)$ and
$$(h, \id_M) \prec \St(\V_0,\St(\U_1)) \prec \St(\V_0,\V_0) = \St(\V_0) \prec \U_0,$$
(1)(ii) implies $h \in U_0$.

$(\ast)_n \, \Rightarrow \, (\ast)_{n+1}$ :
Suppose $h \in \HH(M; K_{n+1})$ satisfies the condition in $(\ast)_{n+1}$.
Since
$$(h|_{K_n}, \id_{K_n}) \prec \St(\V_{n+1},\St(\U_{n+2})) \prec \St(\V_{n+1},\V_{n+1}) =\St(\V_{n+1}) ,$$
the homeomorphism $g=s_{n+1}(h|_{K_n})$ is well-defined and satisfies $(g, \id_M) \prec \U_{n+1}$ by (2)(ii) and $g\in U_{n+1}$ by (1)(ii).

Next we see that $f := g^{-1}h$ satisfies the condition in $(\ast)_n$.
Since $h|_{K_n} = g|_{K_n}$, we have $f := g^{-1}h \in \HH(M; K_n)$.
For each $m \le n$, since
$$(h|_{K_{m-1}} , \id_{K_{m-1}}) \prec \St(\V_m,\St(\U_{n+2})) \mbox{ \ and \ } (g^{-1}, \id_M) \prec \U_{n+1},$$
it follows that
$$(f|_{K_{m-1}} , \id_{K_{m-1}}) \prec \St(\St(\V_m,\St(\U_{n+2})), \U_{n+1})
\prec \St(\St(\V_m, \U_{n+1}), \U_{n+1}) \prec \St(\V_m,\St(\U_{n+1})).$$
The assumption $(\ast)_n$ now implies
$f \in U_n \cdots U_0$ and so $h = gf \in U_{n+1}U_n \cdots U_0$ as required.

Finally, we choose a cover $\V \in \cov(M)$ such that
$$\big\{\St(x,\V):x\in\check M_{n-1}\big\}\prec\V_n\mbox{ \ for all $n\in\w$.}$$
Then, $B(\id_M,\V)=\{h\in\HH_c(M;K):(h,\id_M)\prec\V\}$ is a neighborhood of $\id_M$ in $\HH_c(M;K)$ and
we have $B(\id_M,\V) \subset \RP_{n\in\w}U_n$.
Indeed, any $h \in B(\id_M,\V)$ lies some $\HH(M;K_{n})$ and
the choice of the cover $\V$ implies that
$$(h|_{K_{m-1}},\id_{K_{m-1}})\prec\V_m\prec\St(\V_m,\St(\U_{n+1})) \mbox{ \ for all $m\le n$.}$$
Therefore, the assertion $(\ast)_n$ yields $h \in U_n \cdots U_0$.
\end{proof}

\section{Mapping class groups of non-compact surfaces}\label{sMCG}

In this section,
 we shall prove Theorem~\ref{MCG} announced in Introduction.
For a non-compact connected surface $M$ we need to check the equivalence of the following three conditions:
\begin{enumerate}
\item the mapping class group $\M_c(M)$ is trivial;
\item $\M_c(M)$ is a torsion group;
\item the surface $M$ is exceptional.
\end{enumerate}

Let us recall that a non-compact surface $M$ is exceptional if it is homeomorphic to $(\approx)$ $X\setminus K$ where $X$ is the annulus $\IA$, the disk $\ID$ or the M\"obius band $\IM$ and $K$ is a non-empty compact subset of a boundary circle of $X$.  A subsurface of $M$ means
a subpolyhedron $N$ of $M$ such that
 $N$ is a 2-manifold and $\bd_M N$ is transverse to $\partial M$
 so that $\bd_M N$ is a proper 1-submanifold of $M$ and
 $M \setminus \tint_M N$ is also a 2-manifold.
\smallskip

We shall establish the implications $(1) \Ra (2) \Ra (3) \Ra (1)$,
 the first of which is trivial.
First we verify the implication $(3) \Ra (1)$.

\begin{lemma}
If $M$ is exceptional, then the group $\M_c(M)$ is trivial.
\end{lemma}

\begin{proof}
Take any $h \in \HH_c(M)$. We have to show that $h \in \HH_0(M)$.
In the case $X = \IA$, let $C_0$ denote the boundary circle of $X$ which does not meet $K$.
We can find a compact subsurface $N\subset M$ such that $\supp(h)\subset N$,
$N \approx X$ and $C_0 \subset N$ if $X = \IA$.
Let $C$ denote the boundary circle  of $N$ different from $C_0$.
The frontier $\bd_M N$ of $N$ in $M$ is a compact 1-submanifold of $C$,
which is nonempty by the connectedness of $M$.

Since $N\approx\ID, \IM$ or $\IA$ and $C$ is a boundary circle of $N$,
it follows that $\HH(N,C)$ is path-connected
 (apply \cite[Theorem 3.4]{Ep} in the case $N \approx \IM$).
As is easily observed,
 any $h_1 \in \HH(N,\bd_M N)$ is isotopic to
 some $h_0 \in\HH(N,C)$ in $\HH(N,\bd_M N)$.
Thus, $\HH(N,\bd_M N)$ is also path-connected.
Since $\HH(M,M \setminus \tint_M N) \approx \HH(N,\bd_M N)$,
we have $h \in \HH(M,M\setminus\tint_M N)\subset \HH_0(M)$.
This completes the proof.
\end{proof}

The proof of the implication $(2) \Ra (3)$ is a bit longer
 and is preceded by three lemmas.
The first of them is taken from \cite[Fact 3.2(ii)]{Yag2}.

\begin{lemma}\label{l2}
Every boundary circle $C$ of a non-compact connected surface $M$
 is a retract of $M$.
\end{lemma}

We shall say that an isotopy $h:M\times\II\to M$ has {\em compact support}
 if there is a compact set $K\subset M$ such that
 $h(x,t)=x$ for all $x\in M\setminus K$ and $t\in\II = [0,1]$.

\begin{lemma}\label{l_cpt-supp}
A homeomorphism $h\in\HH(M)$ belongs to $\HH_0(M)$ if and only if
 there is an isotopy $H:M\times\II\to M$ from $h$ to $\id_M$ with compact support.
\end{lemma}

\begin{proof}
The ``if'' part is obvious because
 each isotopy from $h\in\HH(M)$ to $\id_M$ with compact support
 induces a continuous path in $\HH(M)$ from $h$ to $\id_M$,
 which yields $h\in\HH_0(M)$.

Conversely, each $h\in\HH_0(M)$ can be joined to $\id_M$
 by a path $\gamma:\II\to\HH_0(M)$
 (because $\HH_0(M) \approx \IR^\infty\times l_2$ is path-connected).
Since $\gamma(\II)$ is a compact subset of $\HH_c(M)$,
by \cite[Proposition~3.3]{BMSY1}
 $\gamma(\II)\subset\HH(M,M \setminus K)$ for some compact subset $K$ of $M$.
This means that the isotopy induced from $\gamma$ has compact support.
\end{proof}

Our third lemma detects some surfaces $M$ with non-torsion mapping class groups.

\begin{lemma}\label{mcg:l4}
$\M_c(M)$ contains an infinite cyclic subgroup in any of the following cases:
\begin{enumerate}
\item $M$ contains a handle;
\item $M$ contains at least two disjoint M\"obius bands
 (i.e., $M$ contains a Klein bottle with a hole);
\item $M$ contains at least two boundary circles;
\item $M$ contains a M\"obius band and a boundary circle;
\item $M$ is separated by a circle $C \subset {\rm Int}\,M$ into
two non-compact connected subsurfaces $L_1$ and $L_2$.
\end{enumerate}
\end{lemma}

\begin{proof}
(1) Assume that $M$ contains a handle $H$.
For each $n \in \IZ$,
 let $h_n$ denote the $n$-fold Dehn twist
 along the meridian $m$ of the handle $H$.
Then, $h_n \in \HH_0(M)$ iff $n = 0$.
In fact, if $h_n$ is isotopic to $\id_M$,
 then $h_n\ell \simeq \ell$ in $M$ for the longitude $\ell$ of $H$.
Since there is a retraction of $M$ onto $H$ by Lemma \ref{l2},
 we have $h_n\ell \simeq \ell$ in $H$, and this implies that $n = 0$.
\smallskip

(2) When $M$ contains two disjoint M\"obius bands,
 we can connect them by a band to obtain a Klein bottle $K$ with a hole.
We can find two circles $m$ and $\ell$ (a meridian and a longitude) in $K$
 which meet transversely at one point $a$ and
 such that $K$ is obtained by attaching a 2-disk to the wedge $\ell \cup m$
 along the loop $\ell m \ell^{-1}m$ and removing a small disk
 from this disk (cf.\ \cite[p.\,107]{Scott}).
For each $n \in \IZ$,
 let $h_n$ denote the $n$-fold Dehn twist along $m$.
Then, $h_n \in \HH_0(M)$ iff $n = 0$.
In fact, if $h_n$ is isotopic to $\id_M$,
 then $m^n\ell \simeq h_n\ell \simeq \ell$ in $M$.
Since $M$ retracts onto $K$ by Lemma \ref{l2}
 and $K$ retracts onto $\ell \cup m$,
 we have $m^n\ell \simeq \ell$ in $\ell \cup m$.
By contracting $\ell$ onto the point $a$,
 we have $m^n \simeq \ast$ in $m$,
 which implies that $n = 0$.
\smallskip

(3) If $M$ contains two boundary circles $C_1,C_2$,
 then we can connect them by an arc
 and take a regular neighborhood of their union
 to obtain a disk $D$ with two holes,
 which has the two boundary circles $C_1$, $C_2$
 and another boundary circle $C_0 = \bd_M D \subset M\setminus\partial M$.
For each $n \in \IZ$,
 let $h_n$ denote the $n$-fold Dehn twist along a circle in $D$
 parallel to $C_0$.
Then, $h_n \in \HH_0(M)$ iff $n = 0$.
In fact, if $h_n \in \HH_0(M)$, then by Lemma~\ref{l_cpt-supp}
there is an isotopy $h_n \simeq \id_M$ with a compact support $K$.
Take a path $\ell$ in $M$ connecting a point $a \in C_1$
with a point $b \in M \setminus (D \cup K)$ and crossing the circle $C_0$ once.
Then, we have $\ell_n \equiv h_n\ell \simeq \ell$ in $M$,
 where the homotopy keeps the end point $a$ in $C_1$
 and fixes the end point $b$.
Thus, if we retract $M$ onto $D$ and collapse the circle $C_1$ to a point,
 then we obtain an annulus $A$ with the boundary circles $C_2$ and $C_0$
 and a homotopy $\hat{\ell}_n \simeq \hat{\ell}$ rel.\ endpoints in $A$.
This induces a homotopy of the loops $C_0^n \simeq \hat{\ell}_n\hat{\ell}^{-1} \simeq \ast$ in $A$,
which means that $n = 0$.
\smallskip

(4) If $M$ contains a boundary circle and a M\"obius band,
 then we can connected them by an arc
 and take a regular neighborhood of their union
 to obtain a M\"obius band $N$ with a hole.
Then, $C_0 = \bd_M N$ is a circle in $M \setminus \partial M$.
For each $n \in \IZ$,
 let $h_n$ denote the $n$-fold Dehn twist along a circle in $N$
 parallel to $C_0$.
Then, $h_n \in \HH_0(M)$ iff $n = 0$.
Indeed, if $h_n \in \HH_0(M)$, then we can proceed in the same way as in (3)
 to obtain a M\"obius band with the boundary circle $C_0$ and
 a homotopy $C_0^n \simeq \ast$, which means $n = 0$.
 \smallskip

(5)
For every $n \in \IZ$, denote by $h_n$ the $n$-fold Dehn twist along the circle $C$.
Then, $h_n \in \HH_0(M)$ iff $n = 0$.
In fact, if $h_n \in \HH_0(M)$, then there exists an isotopy $h_n \simeq \id_M$ with a compact support $K$.
We can find a path $\ell$ in $M$ connecting a point in $L_1\setminus K$
 with a point of $L_2\setminus K$ and crossing $C$ once.
Then, $h_n\ell \simeq \ell$ in $M$ rel.\ end points and
so $C^n \simeq (h_n\ell) \ell^{-1} \simeq \ast$ in $M$.
Since $M$ retracts onto $C$ by Lemma \ref{l2},
 we have $C^n \simeq \ast$ in $C$, which implies that $n = 0$.
\end{proof}

With Lemma~\ref{mcg:l4} in hands,
 we are now able to prove the implication $(2) \Ra (3)$ of Theorem~\ref{MCG}.

\begin{lemma}
If $\M_c(M)$ is a torsion group,
 then the surface $M$ is exceptional.
 \end{lemma}

\begin{proof}
 Assuming that $\M_c(M)$ is a torsion group and applying Lemma~\ref{mcg:l4},
 we conclude that $M$ contains at most one M\"obius band,
 at most one boundary circle (but not simultaneously),
 no handle, and no circle
 separating $M$ into two non-compact connected subsurfaces.

Write $M$ as the union $M=\bigcup_{n\in\w}M_n$
 of compact connected subsurfaces $M_n$ such that
\begin{enumerate}
\item $M_n \subset \tint_M M_{n+1}$,
\item if $L$ is a connected component of $\check M_n =M \setminus \tint_M M_n$,
 then $L$ is non-compact and $L \cap M_{n+1}$ is connected,
\item if $\partial M \neq \emptyset$ then $M_0 \cap \partial M\neq\emptyset$.
\end{enumerate}

First we shall show the following claim:
\begin{claim}
Every $M_n$ has exactly one boundary circle
 meeting $\check M_n$.
\end{claim}

Each connected component $L$ of $\check M_n$
 meets exactly one boundary circle of $M_n$.
In fact, if $L$ meets two boundary circles $C_1$ and $C_2$ of $M_n$,
 then we can find a circle $C_0$ which meets $C_1$ exactly at one point.
This implies that $M$ contains a handle or a Klein bottle with a hole.
But this contradicts Lemma~\ref{mcg:l4} (1), (2).

Now, assume that $M_n$ has two boundary circles $C_1, C_2$
 meeting $\check M_n$.
Then, the center circle $C$ of a collar neighborhood $A$ of $C_1$ in $M_n$
separates $M$ into two non-compact connected subsurfaces $L_1$ and $L_2$,
 where $L_1$ contains all connected components of $\check M_n$
 meeting $C_1$ and $L_2$ contains $M_n \setminus A$ and
 all remaining connected components of $\check M_n$.
This contradicts Lemma~\ref{mcg:l4} (5).

\medskip
Now, we consider three possible cases.

\smallskip
Case (i): $M$ contains no boundary circle and no M\"obius band.
In this case, each subsurface $M_n$ must be a disk.
In fact, the subsurface $M_n$ contains no handle and no M\"obius band
 because so does $M$.
Then, $M_n$ is homeomorphic to the sphere with $k$ holes.
The connectedness and non-compactness of $M$ imply that $k \geq 1$.
Taking into account Claim and the fact that $M$ has no boundary circle,
 we conclude that $k=1$, which means that each $M_n$ is a disk.

If $\partial M=\emptyset$,
 then $M_n\subset {\rm Int}\, M_{n+1}$ and
 $M_{n+1}\setminus {\rm Int}\, M_n$ is an annulus,
 and hence $M=\bigcup_{n\in\w}M_n$ is homeomorphic to an open disk
 $\ID\setminus \partial\ID$.
If $\partial M\ne \emptyset$,
 then $M_{n+1}\setminus \Int_M  M_n$ consists of disks,
 each of which is attached to $M_n$ along an arc in $\partial M_n$.
This information can be used to show that
 $M$ is homeomorphic to $\ID \setminus K$
 for some non-empty (0-dimensional) compact subset $K \subset \partial\ID$.
\smallskip

Case (ii): $M$ contains  a M\"obius band.
In this case, we may assume that $M_0$ contains a M\"obius band.
Since $M$ contains no boundary circle,
the same argument as in Case (i) implies that every $M_n$ is a M\"obius band.

If $\partial M = \emptyset$,
 then each $M_{n+1}\setminus {\rm Int}\, M_n$ is an annulus
 and then $M$ is homeomorphic to an open M\"obius band
 $\IM \setminus \partial\IM$.
If $\partial M \neq \emptyset$,
 then $M_{n+1}\setminus\tint_M M_n$ consists of disks,
 each of which is attached to $M_n$ along an arc in $\partial M_n$.
In this case, $M$ is homeomorphic to $\IM\setminus K$,
 where $K$ is a non-empty compact subset of $\partial\IM$.
\smallskip

Case (iii): $M$ contains a boundary circle $C$.
In this case, we may assume that $C$ belongs to $M_0$.
Since $M$ contains no M\"obius band,
we conclude that every $M_n$ is an annulus.
Let $C_1$ and $C_2$ denote the boundary circles of the annulus $\IA$.

If $\partial M = C$,
 then each $M_{n+1}\setminus \tint_M M_n$ is an annulus
 and the pair $(M,C)$ is homeomorphic to $(\IA \setminus C_2,C_1)$.
If $\partial M \neq C$,
 then we may assume that $M_0$ meets $\partial M \setminus C$.
In this case, $M_{n+1} \setminus \tint_M M_n$ consists of disks,
 each of which is attached to $M_n$ along an arc in $\partial M_n \setminus C$, and hence
 the pair $(M,C)$ is homeomorphic to $(\IA \setminus K,C_1)$,
 where $K$ is a non-empty compact subset of $C_2$.
\end{proof}

\section*{Acknowledgment}

The authors express their sincere thanks to Sergiy Maksymenko and B\l a\.zej  Szepietowski,
 for fruitful and stimulating discussions
 related to  mapping class groups.

\end{document}